\documentclass[a4paper,11pt,fleqn]{amsart}
\usepackage[utf8]{inputenc}
\usepackage{epsfig,graphics}
\usepackage[all,knot,poly,color,curve]{xy}
\usepackage{array}
\usepackage{amsthm,amssymb,amsbsy,bbm,a4wide,verbatim,url}
\usepackage[marginpar]{todo}
\usepackage{longtable}
\usepackage{rotating}
\usepackage{color,fdsymbol}

\usepackage{fancybox}
\usepackage{enumitem}

\def\un{\mathbbm{1}}

\def\PSp{\mathrm{PSp}}

\def\SU{\mathrm{SU}}

\def\diag{\mathrm{diag}}

\def\F{\mathbbm{F}}

\def\sl{\mathfrak{sl}}

\def\Ker{\mathrm{Ker}}

\def\GL{\mathrm{GL}}

\def\SL{\mathrm{SL}}

\newtheorem{theorem}{Theorem}[section]

\numberwithin{equation}{section}

\date{February 9, 2020}

\author{Ivan Marin}
\title{On the largest representation in type $H_4$}
\dedicatory{To the memory of Kay Magaard}

\begin{document}

\maketitle

\begin{abstract}
In this technical note, we complete the PhD work of A. Esterle about determining the image of any Artin group of finite Coxeter type inside the associated Hecke algebra over a finite field, when the latter is semisimple. The only remaining case was the 48-dimensional irreducible
representation in type $H_4$, for which the image is proven here to be $\Omega_{48}^+$.
\end{abstract}

\section{Introduction}

In \cite{THESTERLE}, A. Esterle determined the image of the Artin groups of finite Coxeter type
inside the corresponding Hecke algebra over a finite field, under some conditions on the
parameter essentially made so that the latter is semisimple. There is just \emph{one} case for which
the methods of \cite{THESTERLE} where not conclusive, because a main technical tool due
to Guralnick and Saxl could not be applied there. This case is the largest irreducible representation
of the exceptional Coxeter type $H_4$. Here we solve this remaining case by going back to the
original basic tool in this area, namely Aschbacher's classification theorem.

We use the notational conventions of A. Esterle's thesis \cite{THESTERLE}. In particular, as in  \cite{THESTERLE}, let us denote $A_{H_4}$ the Artin group of type $H_4$ and $\mathcal{A}_{H_4}$
its derived subgroup. We consider the representation $\rho : \mathcal{A}_{H_4} \to \GL_{48}(\overline{\F_p})$
deduced from the only 48-dimensional irreducible representation of the Hecke
algebra in type $H_4$, and we want to determine the image $G$ of $\rho$. The convention on Hecke algebras is that the Artin generators have eigenvalues $-1,\alpha$.

It is the only unsolved case in \cite{THESTERLE}, where it is shown that $G$ is
a subgroup of $\Omega_{48}^+(\sqrt{q})$, depending on the parameter
$\alpha \in \F_q^{\times}$ with $\F_q = \F_p(\alpha)$. It is assumed that  $p \not\in \{2,3,5 \}$ and the order of $\alpha$ in $\F_q^{\times}$ is assumed not to divide 20,30 or 48. These conditions imply $q \geq 19$ and, when
$q$ is a square, that $q \geq 121$.
Moreover, since the other cases are already dealt with in \cite{THESTERLE}, one can assume that (condition '$1 \sim 2$' in \cite{THESTERLE}) 
there exists a field automorphism $\Phi$ of $\F_p(\alpha,\xi + \xi^{-1})$
such that $\Phi(\alpha + \alpha^{-1}) = \alpha + \alpha^{-1}$ and
$\Phi(\xi + \xi^{-1}) = \xi^2 + \xi^{-2})$, which is equivalent (\cite{THESTERLE} lemma 8.1.2) to the condition $\F_p(\alpha,\xi+\xi^{-1}) = \F_p(\alpha+\alpha^{-1},\xi+\xi^{-1})$. We set $\F_g = \F_p(\alpha+\alpha^{-1})$.

Because of the existence of $\Phi \in \mathrm{Gal}(\F_g(\xi + \xi^{-1})/\F_g)\setminus \{ \mathrm{Id} \}$, $\F_g$ has index 2 inside $\F_g(\xi + \xi^{-1})$. Now, it is proven in \cite{THESTERLE} that the representation of $\mathcal{A}_{H_4}$ under consideration admits a matrix model over $\F_g(\xi + \xi^{-1})$, so we can write it as $\rho : \mathcal{A}_{H_4} \to \GL_{48}(\F_g(\xi + \xi^{-1}))$. If $\F_p(\alpha) = \F_g$, that is $q = g$, since there is only
one 48-dimensional irreducible representation of $A_{H_4}$ that factorizes through the Hecke algebra of type $H_4$, and since this one is defined over $\F_g$,  
we have $\Phi \circ \rho \simeq \rho$ and we can assume that $\rho$ takes value inside $\GL_{48}(\F_g) = \GL_{48}(\F_q)$ (see e.g. \cite{THESTERLE} lemma 3.2.5). If not,
we have that $[\F_p(\alpha): \F_g] = 2$ hence $\F_p(\alpha) = \F_g(\xi + \xi^{-1})$ and $\rho$ also takes value inside $\GL_{48}(\F_q)$, this time with
$q = g^2$. But then $\Phi$ is the only non-trivial element of $\mathrm{Gal}(\F_p(\alpha)/\F_g)$, hence exchanges $\alpha$ and $\alpha^{-1}$.
From Proposition 2.1.2 of \cite{THESTERLE} this implies that $\Phi \circ \rho \simeq \rho$ hence $\rho$ can actually be defined over $\F_p(\alpha+\alpha^{-1}) = \F_g$ and this proves that $\rho$ can be defined with values in $\GL_{48}(\F_g)$ again.

The purpose of this note is to establish the following statement, which solves conjecture 8.5.1 of \cite{THESTERLE} and completes the goal of A. Esterle thesis. Notice that actually the case $\F_p(\alpha) = \F_p(\alpha+\alpha^{-1})$ was overlooked there, which makes Conjecture 8.5.1 trivially false as stated in this case,
so that Theorem 8.5.1 of \cite{THESTERLE} needs to be corrected. Explicitely, in the cases 1.a, 1.d of the theorem, the factor $\Omega_{48}^+(\sqrt{q})$ of the decomposition must be replaced by $\Omega_{48}^+(q)$.

\begin{theorem} In the case $1 \sim 2$, and with the notations of \cite{THESTERLE}, we have $G = \Omega_{48}^+(g)$, with 
$\F_g = \F_p(\alpha+\alpha^{-1})$. 
\end{theorem}

Distinguishing cases, this statement can be rewritten as follows.
\begin{enumerate}
\item If $\F_q = \F_p(\alpha) = \F_p(\alpha+\alpha^{-1})$, then $G = \Omega_{48}^+(q)$
\item If $\F_q = \F_p(\alpha) \neq \F_p(\alpha+\alpha^{-1}) = \F_{\sqrt{q}}$, then $G = \Omega_{48}^+(\sqrt{q})$
\end{enumerate}

As explained in \cite{THESTERLE}, this case could not be tackled by the
methods used there, because the restrictions to some parabolic subgroup do
not contain elements of the type needed for applying a very handy theorem of Guralnick and Saxl. Therefore
we change strategy and directly apply Aschbacher's theorem (see \cite{ASCHBACHER}). 
By the considerations above this theorem implies that it is sufficient
to show that $G$ cannot belong to any of the Aschbacher classes $\mathcal{C}_1,\dots,\mathcal{C}_7$ and $\mathcal{S}$.

By \cite{THESTERLE} theorem 8.1.1 p. 166 and the
fact that the representation $48_{rr}$ of the Coxeter group of type $H_4$ restricts
to $3_s + 3'_s + \overline{3_s} +  \overline{3'_s} + 2 \times 4_r + 2 \times 4'_r + 2 \times 5_r + 2 \times 5'_r$
as a representation of type $H_3$ (see \cite{THESTERLE} Table 8.3),
we know that $G$ contains
$Q_3 \times Q_4 \times Q_5$ with  $Q_4 \in \{ \SL_4(q), \SU_4(\sqrt{q}) \}$
and $(Q_3,Q_5) \in \{ (\SL_3(q^2), \SL_5(q)),(\SL_3(q),\SU_5(\sqrt{q})) \}$ embedded inside $\SL_{48}(\F_q(\xi+\xi^{-1}))$ via
$$
(q_3,q_4,q_5) \mapsto (q_3,\bar{q}_3, ^t q_3^{-1}, ^t\bar{q}_3^{-1},q_4,q_4, ^t q_4^{-1}, ^t q_4^{-1}, q_5,
q_5, ^t q_5^{-1}, ^t q_5^{-1})
$$
where $\bar{x} = \Phi(x)$. We also know that this subgroup normally generates $G$ (\cite{THESTERLE} Lemma 8.2.1, p. 170).
We have 
$|Q_3| \geq \max(|\SL_3(19)|,|\SL_3(19^2)|) \geq 1.5\times 10^{19}$, 
$|Q_4| \geq \max(|\SL_4(19)|, |\SU_4(11)|) \geq 4\times 10^{15}$, 
$|Q_5| \geq \max(|\SL_5(19)|,|\SU_5(11)|) \geq 9.7\times 10^{24}$. 
We also notice that the order of the $p$-Sylow subgroups of $Q_3,Q_4,Q_5$ are at least $q^3$,$\sqrt{q}^6$ or $q^6$,
and $\sqrt{q}^{10}$ or $q^{10}$, respectively.

When one needs to distinguish the two cases, we excerpt again from \cite{THESTERLE} Theorem 8.1.1 that,
when $\F_q = \F_p(\alpha) = \F_p(\alpha+\alpha^{-1})$, then $(Q_ 3,Q_5)= ( \SL_3(q^2), \SL_5(q))$
and when $\F_q = \F_p(\alpha) \neq \F_p(\alpha+\alpha^{-1}) = \F_{\sqrt{q}}$, then
$(Q_ 3,Q_5) = ( \SL_3(q), \SU_5(\sqrt{q}))$.

\medskip

{\bf Acknowledgements.} I thank A. Esterle for numerous discussions, O. Brunat for a careful reading, F. L\"ubeck, G. Hiss and G. Malle for references.

\section{Classical Aschbacher classes}
Assume $G$ is contained in a maximal subgroup $\Gamma$ of $\GL_{48}(r)$ belonging to
Aschbacher's class $\mathcal{C}_i$ with $1 \leq i \leq 7$ as in \cite{ASCHBACHER}.
Class $\mathcal{C}_1$ is excluded because $G$ acts absolutely irreducibly (\cite{THESTERLE}, proposition 8.2.2).
Similarly, classes $\mathcal{C}_6$ and $\mathcal{C}_7$ are excluded because $48$ is not a non-trivial power.
We thus need to consider classes $\mathcal{C}_2,\mathcal{C}_3,\mathcal{C}_4,\mathcal{C}_5$.
\subsection{Class $\mathcal{C}_2$}

We assume now that $\Gamma$ lies in class $\mathcal{C}_2$, that is $\F_r^{48} = V = V_1\oplus \dots \oplus V_s$
with $\dim V_1 = \dots = \dim V_s = v$ so that $48 = vs$, and $\Gamma = \GL_v(r) \wr \mathfrak{S}_s$
with $v,s \geq 2$.
Let $\pi : \Gamma \to \mathfrak{S}_s$ denote the associated projection and consider
the map $\pi_k : Q_k \to \mathfrak{S}_s$ obtained by composing $\pi$ with $Q_k \to Q_3 \times Q_4 \times Q_5 \to G$. Since $Q_k$ is almost simple, either $\Ker \pi_k < Z(Q_k)$ or $\Ker\pi_k = Q_k$. In the
first case, $\pi_k$ induces an embedding of the $p$-Sylow subgroup of $Q_k$ inside $\mathfrak{S}_s$. But
since $p \geq 11$ and $s$ is a proper divisor of $48$, an immediate computation shows that the orders of the
$p$-Sylow subgroups of $\mathfrak{S}_s$ for such $s$ are at most $p^2$, contradicting our conditions. 
It follows that $\pi(G) = \{ 1\}$ hence $G$ acts non-irreducibly on $V$, a contradiction.

\subsection{Class $\mathcal{C}_3$} We have $\Gamma < \SL_{48}(\F_g)$. If $\Gamma$ 
belongs to the class $\mathcal{C}_3$, then we can write $\F_g^{48} = \F_s^{h}$
for $h < 48$, with $\F_s \subsetneq \F_g$,
and $\Gamma < \SL_{h}(\F_s) . C_m$ for
some cyclic group $C_m$.
Since $Q_3 \times Q_4 \times Q_5$ is perfect, its
image 
inside $C_m$ is trivial, whence its image inside $G< \Gamma$ lies inside $\SL_{h}(\F_s)$. Since
it normally generates $G$, we get that $G < \SL_{h}(\F_s)$. Since the 48-dimensional representation of $G$ is
absolutely irreducible 
we get a contradiction.

\subsection{Class $\mathcal{C}_4$}
We want to prove that the action of $G$
is tensor-indecomposable. For this it is sufficient to prove
that the action of $Q_3< G$ is tensor-indecomposable. We
thus consider $W = V \otimes_{\F_r} \overline{\F}_p  = \overline{\F}_p^{48}$
as a $Q_3$-module, and assume it can be written as $W_1 \otimes W_2$.
It is a semisimple module, that can be decomposed as $U + \overline{U} + U^* + \overline{U}^* + 36.\un$
where $U,\overline{U},U^*,\overline{U}^*$ are four absolutely irreducible pairwise non-isomorphic
3-dimensional representations of $Q_3$. Since the only prime divisors of $48$ are $2$ and $3$
we have $48 \not\equiv 0 \mod p$ hence $W_1$ and $W_2$ are semisimple by \cite{SERRE} Theorem 2.4.
We write a direct sum decomposition
$W_k = m_k \un + \sum_{i=1}^{r_k} a_i^{(k)} F_i^{(k)}$ where the $F_i^{(k)}$ are pairwise
nonisomorphic nontrivial simple modules. Since $Q_3$ is perfect they all have dimension at least $2$. Then, the fact that the nontrivial irreductible constituents of $W$
appear multiplicity free implies $m_1,m_2 \in \{ 0, 1 \}$. We have
$$
W = W_1 \otimes W_2 = m_1m_2 \un + m_1 \sum_{i=1}^{r_2}a_i^{(2)} F_i^{(2)}
+ m_2 \sum_{i=1}^{r_1}a_i^{(1)} F_i^{(1)}
+   \sum_{i=1}^{r_1}\sum_{k=1}^{r_2}a_i^{(1)}a_j^{(2)} F_i^{(1)} \otimes F_i^{(2)}
$$
and, since $W$ contains $36.\un$ and $F_i^{(1)} \otimes F_i^{(2)}$ may contain the
trivial representation at most once by Schur's lemma, this implies
$$
m_1 m_2 + \left( \sum_{i=1}^{r_1}a_i^{(1)}\right)
\left( \sum_{i=1}^{r_2}a_i^{(2)}\right) \geqslant 36
$$
and $m_1m_2 \leq 1$ implies 
$\left( \sum_{i=1}^{r_1}a_i^{(1)}\right)
\left( \sum_{i=1}^{r_2}a_i^{(2)}\right) \geq 35$.
We have $\dim F_i^{(k)} \geq 2$ for all $i,k$ since $Q_3$ is perfect hence $\un$ is the only 1-dimensional representation. Then $48 = (\dim W_1)(\dim W_2)$
and $\dim W_k = 
m_k  + \sum_{i=1}^{r_k} a_i^{(k)} \dim F_i^{(k)} \geq
   2 \sum_{i=1}^{r_k} a_i^{(k)}$ thus $48 = \dim W \geq 4 \times 35$, a contradiction. Therefore
it is tensor-indecomposable.

\subsection{Class $\mathcal{C}_5$}

We have $\Gamma < \SL_{48}(\F_g)$. If $\Gamma$ 
belongs to the class $\mathcal{C}_5$, then we can write $\F_g^{48} = \F_s^{48} \otimes_{\F_s} \F$ with $\F_s \subsetneq \F$, and $\Gamma < D = \SL_{48}(\F_s).(\F_g)^{\times}$.
Since $D/\SL_{48}(\F_s) \simeq (\F_g)^{\times})/(\F_{s})^{\times}$ is a cyclic
group and $Q_3\times Q_4 \times Q_5$  is perfect, its
image 
inside $D/\SL_{48}(\F_s)$ is trivial, whence its image inside $G< \Gamma$ lies inside $\SL_{48}(\F_s)$. Since
it normally generates $G$, we get that $G < \SL_{48}(\F_s)$.

We then distinguish two cases. First assume that $\F_q = \F_p(\alpha) \neq \F_p(\alpha+\alpha^{-1}) = \F_{\sqrt{q}} = \F_g$. In that case, from \cite{THESTERLE} theorem 8.1.1 (p.166) we get
that $Q_3 = \SL_3(q)$. Let $q_3 = \diag(\alpha,\alpha^{-1}, 1) \in Q_3$, $q_4 = 1$, $q_5 = 1$. The image of $(q_3,q_4,q_5)$ inside $G$ has trace $2(\alpha + \alpha^{-1})
+ 2 \Phi(\alpha + \alpha^{-1})+40= 4 (\alpha + \alpha^{-1})+40$.
It follows that $\F_s \supset \F_p(\alpha + \alpha^{-1}) = \F_{\sqrt{q}}$
hence $s = \sqrt{q}$, contradicting $\F_s \subsetneq \F_{\sqrt{q}}$.

We then assume $\F_q = \F_p(\alpha) = \F_p(\alpha+\alpha^{-1})= \F_g$. In that case, from \cite{THESTERLE} theorem 8.1.1 (p.166) we get
that $Q_3 = \SL_3(q^2)$ and $Q_5 = \SL_5(q)$. Letting
$q_5 = \diag(\alpha,\alpha^{-1},1,1,1)$ and $q_3 = 1$, $q_4 = 1$, we get that the image of $(q_3,q_4,q_5)$ inside $G$ has trace $4(\alpha+\alpha^{-1}) + 40$ hence $\F_s \supset \F_p(\alpha + \alpha^{-1}) = \F_q$,
contradicting $\F_s \subsetneq \F_{q}$.

\section{Aschbacher class $\mathcal{S}$}

We finally exclude the groups in the $\mathcal{S}$ class.

\subsection{Groups of Lie type in non-defining characteristic and sporadic groups}
We use \cite{HISSMALLE} in this part where all the cases in small dimension
have been classified.

$$
\begin{array}{l|l|r|l||c|}
\hline
& \Gamma & \ell & \mbox{field} & |\Gamma| \\
\hline
(1) & 2.\mathfrak{A}_8 & \neq 2 & & 40320\\
(2) & \mathfrak{A}_9 & 0,2& & 181440 \\
(3) & 2.\mathfrak{A}_9 & 3& & 362880 \\
(4) & 2.\mathfrak{A}_9 & \neq 2,3 & i6 &362880
\\
(5) & \mathfrak{A}_{10} & 2 & & 1814400\\
(6) & 2.\mathfrak{A}_{10} & 3 & &  3628800 \\
(7) & 2.\mathfrak{A}_{10} & \neq 2,3 & i6&3628800 \\
(8) & 12_1.L_3(4) & 0,7 & z12,b5 & 241920 \\
(9) & 12_2.L_3(4) & 0,7 & z12,b5 & 241920 \\
(10) & 12_2.L_3(4) & 5 & z12 & 241920 \\
(11) & 3.U_3(5) & \neq 3,5 & z3 & 378000 \\
(12) & 2.S_6(2) & \neq 2,7 &  & 2903040\\
(13) & O_8^+(2) & 3 &  & 174182400 
\\
(14) & 2.Sz(8) & 5 & c13 & 58240 \\
(15) & 12.M_{22} & 5 & z12,b11 & 5322240 \\
\hline
\end{array}
$$
This table is the relevant excerpt of exceptional cases from table 3 in \cite{HISSMALLE} (corrected as table 2 in \cite{HISSMALLEERR}),
completed by the order of the group.
Since $|\Gamma| \geq |G| \geq |Q_3| \times |Q_4|\times |Q_5| \geq 58\times 10^{58}$, the computation of this order is sufficient to
dismiss these cases.

The generic cases explained in \cite{HISSMALLE} (table 2 there) are of two types. One of these types (cases (b),(c),(d) in \cite{HISSMALLE}, table 2) is when $G$ is either $2.PSL_2(m)$
or $PSL_2(m)$, with $m \leq 2\times 48+1 \leq 100$.
But then 
$$
|G| \leq 2 \times 100 \times (100^2 - 1 ) \leq 2 \times 10^6
$$
contradicting $|G| \geq 58\times 10^{58}$.
The other type (case (a) in \cite{HISSMALLE}, table 2)
is when $G = \mathfrak{A}_n$ with $n \in \{ 49, 50 \}$.
Then, since $p \geq 7$, 
$$v_p(|G|) \leqslant v_p(50!) = \lfloor \frac{50}{p^2} \rfloor +
\lfloor \frac{50}{p} \rfloor \leq 2 + 8 = 10  
$$
while $Q_3 \times Q_4 \times Q_5$ has a $p$-Sylow of order at least $p^{19}$, dismissing this case again.

Therefore this case is excluded, too.

\subsection{Groups of Lie type in natural characteristic}

We use \cite{LUBECK} here. In the tables there (appendices A...)
the representations of dimension $48$ (except for the natural representation of the classical groups) which appear
are the following ones :
\begin{enumerate}
\item case $A_2$, $p \neq 7$
\item case $A_6$, $p \neq 7$
\item case $B_3$, $p \neq 7$
\item case $B_4$, $p =2$
\item case $C_3$, $p =2$
\item case $C_4$, $p \neq 3$
\item case $D_4$, $p =2$
\end{enumerate}
Also, it is not taken into account there the type $A_1$, which 
we can rule out separately, as in this case $G$ would have an abelian $p$-Sylow subgroup which contradicts
the fact that the one of $Q_3 \times Q_4 \times Q_5$ is not.

Since $p \neq 2$, the only cases to consider are then (1),(2),(3),(6). 
In all cases $Q_3,Q_4,Q_5$ have $p$-Sylow subgroups of order $g^6,g^6,g^{10}$. It follows that a $p$-Sylow subgroup of $Q_3 \times Q_4 \times Q_5$ has order $g^{22}$. In cases (1),(2),(3),(6) (and $p \neq 2$) $G$ is contained inside $\SL_3(g)$,$\SL_7(g)$,$\Omega_7(g)$,$\PSp_8(g)$, which
have $p$-Sylow subgroups of order $g^{3}$,$g^{21}$,$g^9$,$g^{16}$, respectively, a contradiction.


\begin{thebibliography}{DWKL}
\bibitem{ASCHBACHER} M. Aschbacher, {\it On the maximal subgroups of the classical groups}, Invent. math. {\bf 76} (1984), 469-514.
\bibitem{THESTERLE} A. Esterle, {\it Artin groups and Hecke algebras over finite fields}, doctoral thesis, Amiens, 2018, arxiv:1808.03687v1.
\bibitem{HISSMALLE} G. Hiss, G. Malle, {\it Low-dimensional representations of quasi-simple groups},
 LMS J. Comput. Math. {\bf 4} (2001), 22-63. 
\bibitem{HISSMALLEERR} G. Hiss, G. Malle, {\it Corrigenda :  Low-dimensional representations of quasi-simple groups},
  LMS J. Comput. Math. {\bf 5} (2002), 95-126.
\bibitem{LUBECK} F. L\"ubeck, {\it Small degree representations of finite Chevalley group in defining characteristic},
 LMS J. Comput. Math. {\bf 4} (2001), 135-169. 
\bibitem{MBM} K. Magaard, O. Brunat, I. Marin, {\it Image of the braid groups inside the finite Iwahori-Hecke algebras}, J. Reine Angew. Math. {\bf 733} (2017), 161-182. 
\bibitem{MARIN} I. Marin, {\it Infinitesimal Hecke Algebras II}, preprint 2009, \verb+arxiv:0911.1879v1+ 
\bibitem{SERRE} J.P. Serre, {\it Semisimplicity and Tensor Products of
Group Representations: Converse Theorems}, J. Algebra {\bf 194}, 496-520 (1997).
\end{thebibliography}
\end{document}